\def \beq{\begin{equation}}
\def \eq{\end{equation}}
\def \berr{\begin{eqnarray}}
\def \err{\end{eqnarray}}
\def \nn{\nonumber}
\def \a{\alpha}
\def \b{\beta}

\def \d{\delta}
\def \r{\rho}

\def \Del{\Delta}
\def \w{\omega}

\def \eps{\varepsilon}
\def \r{\rho}

\def \l{\lambda}

\def \U{{\cal U}}
\def \W{{\cal W}}

\def \Par{\mbox{Par}}
\def \ch{\mbox{ch}}

\def \({\left(}
\def \){\right)}
\def \<{\langle}
\def \>{\rangle}
\def \[{\left[}
\def \]{\right]}
\def \obar{\overline}

\def\tens{\mathop{\otimes}}

\def\uk{\underline{k}}
\def\ttens{\tilde{\tens}}

\newcommand \reals{I \! \! R}
\newcommand \compl{C \! \! \! \! {\scriptscriptstyle {}^{{}_|}}\ }
\newcommand \N{I \! \! N}
\newcommand \R{I \! \! R}
\newcommand \Z{Z \! \! \! Z}

\def\reps{representations }
\def\rep{representation }

\documentstyle[12pt]{article}
\input epsf

\newtheorem{prop}{Proposition}[section]
\newtheorem{theorem}[prop]{Theorem}
\newtheorem{lemma}[prop]{Lemma}
\newtheorem{definition}[prop]{Definition}

\evensidemargin \oddsidemargin
\begin{document}

\begin{titlepage}
\begin{center}
        \hfill  LBL-39526 \\
          \hfill    UCB-PTH-96/46 \\
\vskip .5in

{\large Finite dimensional unitary representations of quantum
Anti--de Sitter groups at roots of unity} 
\vskip .50in

Harold Steinacker   \footnote{email:
hsteinac@physics.berkeley.edu}

\vskip .2in
{\em

Theoretical Physics Group \\
Ernest Orlando Lawrence Berkeley National Laboratory  \\
University of California, Berkeley, California  94720  \\

and

Department of Physics  \\
University of California, Berkeley, California  94720 }

\end{center}

\vskip .5in

\begin{abstract}

We study irreducible unitary \reps of $U_q(SO(2,1))$ 
and $U_q(SO(2,3))$ for $q$ a root of unity, which are 
finite dimensional. Among others, unitary \reps corresponding to all
classical one--particle representations 
with integral weights are found
for $q = e^{i \pi /M}$, with $M$ being large enough. 
In the "massless" case with spin bigger than or equal to
1 in 4 dimensions, they are unitarizable only
after factoring out a subspace of "pure gauges", as classically.
A truncated associative tensor product describing unitary 
many--particle representations is defined for 
$q = e^{i\pi /M}$.

\end{abstract}

\end{titlepage}

\section{Introduction}

In recent years, the development of Noncommutative Geometry 
has sparked  
much interest in formulating physics and in particular quantum field 
theory on quantized, i.e. noncommutative spacetime. The idea 
is that if there are no more "points" in spacetime, such a theory 
should be well behaved in the UV.

Quantum groups \cite{FRT,jimbo,drinfeld}, although discovered in a 
different context, can be understood as generalized "symmetries" of certain 
quantum spaces. Thinking of elementary particles as irreducible unitary 
representations of the Poincar\'e 
group, it is natural to try to formulate a quantum field 
theory based on some quantum Poincar\'e group, i.e. on some quantized
spacetime. 

There have been many attempts (e.g. \cite{poincare,lukierski}) in this 
direction. One of the difficulties with many versions of a 
quantum Poincar\'e group comes from the fact that the classical 
Poincar\'e group is not semisimple. This forbids  using the well
developed theory of (semi)simple quantum groups, which is e.g. reviewed in 
\cite{CH_P,lusztig_book,jantzen}. In this paper, we 
consider instead the quantum Anti--de Sitter group 
$U_q(SO(2,3))$, 
resp. $U_q(SO(2,1))$ in 2 dimensions, thus taking advantage of much 
well--known mathematical machinery. In the classical case, these groups 
(as opposed to e.g. the de Sitter group $SO(4,1)$) are known to have
positive--energy representations for any spin \cite{fronsdal},
and e.g. allow supersymmetric extensions \cite{zumino}. 
Furthermore, one could argue that the usual choice of flat spacetime 
is a singular choice, perhaps subject to some mathematical artefacts.

With this motivation, we study unitary representations 
of $U_q(SO(2,3))$. Classically, all unitary
representations are infinite--dimensional since the group is 
noncompact. It is well known that at roots of unity, the irreducible 
representations (irreps) of quantum groups are finite dimensional. 
In this paper, we determine 
if they are unitarizable, and show in particular that 
for $q = e^{i\pi /M}$, all 
the irreps with positive energy and 
integral weights are unitarizable, as long as the rest energy $E_0 \geq s+1$ 
where $s$ is the spin, and $E_0$ is below some 
($q$ -- dependent, large) limit. 
There is an intrinsic high--energy cutoff, and only finitely many 
such "physical" representations exist for given $q$. 
At low energies and for $q$ close 
enough to 1, the structure is the same as in the classical case.
Furthermore, unitary \reps exist only at roots of unity
(if $q$ is a phase). For generic roots of unity, their weights are 
non--integral. Analogous results are found for $U_q(SO(2,1))$.
In general, there is a cell--like structure of unitary \reps
in weight space.

In the  "massless" case, the  naive \reps 
with spin bigger than or equal to 1
are reducible and contain a null--
subspace corresponding to "pure gauge" states. It is shown that 
they can be consistently factored out to obtain unitary 
representations with only 
the physical degrees of freedom ("helicities"), as in the classical case 
\cite{fronsdal}.

The existence of finite--dimensional unitary representations of
noncompact quantum groups at roots of unity has already been pointed out
in \cite{dobrev}, where several representations of $U_q(SU(2,2))$ 
and $U_q(SO(2,3))$ (with multiplicity of weights equal
to one) are shown to be unitarizable.
In the latter case they correspond to the Dirac singletons
\cite{dirac}, which are recovered here as well.

We also show that the class of "physical" (unitarizable)
representations is closed
under a new kind of associative truncated tensor product
for  $q = e^{i \pi /M}$, i.e. there
exists a natural way to obtain many--particle
representations.

Besides being very encouraging from the point of view
of quantum field theory, this shows again
the markedly different properties of quantum groups at 
roots of unity from the case of generic $q$ and $q=1$. 
The results are clearly
not restricted to the groups considered here and should 
be of interest on purely mathematical grounds as well.
We develop a method to investigate the structure
of representations of quantum groups at roots of unity
and determine the structure of a large class of \reps of 
$U_q(SO(2,3))$. Throughout 
this paper, $U_q(SO(2,3))$ will be equipped with a non--standard 
Hopf algebra star structure.

The idea to find a quantum Poincar\'e group
from $U_q(SO(2,3))$  is not new: Already in \cite{lukierski}, the 
so--called $\kappa$--Poincar\'e group was constructed by a contraction of
$U_q(SO(2,3))$.
This contraction however essentially takes $q 
\rightarrow 1$ (in a nontrivial way) and destroys the properties of the 
representations 
which we emphasize, in particular the finite dimensionality. 

Although it is not considered here, we want to mention that 
there exists a (space of functions on) 
quantum Anti--de Sitter space on which $U_q(SO(2,1))$ resp. $U_q(SO(2,3))$ 
operates, with  an intrinsic 
mass parameter $m^2 = i(q-q^{-1})/R^2$ where $R$ is the "radius" of 
Anti--de Sitter space (and the usual Minkowski signature for $q=1$) 
\cite{thesis}. 

This paper is organized as follows: In section 2, we investigate the 
unitary representations of $U_q(SO(2,1))$, and define a truncated tensor 
product. In section 3, the most important facts about 
quantized universal enveloping algebras of higher rank are reviewed. 
In section 4, we consider $U_q(SO(5))$ and $U_q(SO(2,3))$,
determine the structure of the relevant irreducible \reps 
(which are finite dimensional) and investigate which ones are
unitarizable.
The truncated tensor product is generalized to the case of $U_q(SO(2,3))$. 
Finally we conclude and look at possible further developments.

\section{Unitary representations of $U_q(SO(2,1))$}

We first consider the simplest case of $U_q(SO(2,1))$, 
which is a real 
form of $\U \equiv U_q(Sl(2,\compl))$, the Hopf algebra defined by
\cite{FRT,jimbo}
\berr
[H, X^{\pm}] &=& \pm 2X^{\pm}, \quad [X^+, X^-] = [H] \\
\Del(H) &=& H\tens 1 + 1 \tens H, \nonumber  \\
\Del(X^{\pm}) &=& X^{\pm} \tens q^{H/2} +          
                            q^{-H/2} \tens X^{\pm},   \nonumber \\
S(X^+) &=& -q X^+, \quad S(X^-) = -q^{-1} X^-, \quad S(H) =-H \nonumber \\
\eps(X^{\pm}) &=& \eps(H) =0 \nonumber
\err
where $[n] \equiv [n]_q = \frac{q^n-q^{-n}}{q-q^{-1}}$.
To talk about a real form of $U_q(SL(2,\compl))$, one has to impose a 
reality condition, i.e. a star structure, and there may be several 
possibilities.
Since we want the algebra to be implemented by a unitary 
representation on a Hilbert space, the star operation should be
an antilinear antihomomorphism of the algebra. 
Furthermore,  we will 
see that to get finite dimensional unitary representations, $q$ must 
be a root of unity, so $|q|=1$.
Only at roots of unity the representation 
theory of quantum groups differs essentially from the classical case, 
and new features such as finite dimensional unitary representations 
of noncompact groups can appear. 
This suggests the following
star  structure corresponding to $U_q(SO(2,1))$:
\beq
H^{\ast} = H, \quad (X^+)^{\ast} = -X^-   \label{star_sl2}
\eq
whis is simply 
\beq
x^{\ast} = e^{-i\pi H/2} \theta(x^{c.c.}) e^{i\pi H/2}  \label{invol_help}
\eq
 where $\theta$ is the usual (linear)  Cartan--Weyl involution
and $x^{c.c.}$ is the complex conjugate of $x\in 
\U$. Since $q$ is a phase, $q^{c.c.} = q^{-1}$, and 
\beq
(\Del(x))^{\ast} = \Del(x^{\ast})  \eq
provided
\beq
(a \tens b)^{\ast} = b^{\ast} \tens a^{\ast}.  \label{Del_star}
\eq
Then $(S(x))^{\ast} = S(x^{\ast})$, which is a 
non--standard Hopf algebra star structure. In particular, 
(\ref{Del_star}) is chosen as e.g. in \cite{mack_schom}, which is 
different from the standard definition.
Nevertheless, this is 
perfectly consistent with a many--particle
interpretation in Quantum Mechanics or Quantum Field Theory
as discussed in \cite{thesis}, where it is shown e.g. how
to define an invariant inner product on the tensor product 
with the "correct" classical limit.

\begin{figure}
 \epsfxsize=4in
  \vspace{-2in} 
   \hspace{0.8in}
   \epsfbox{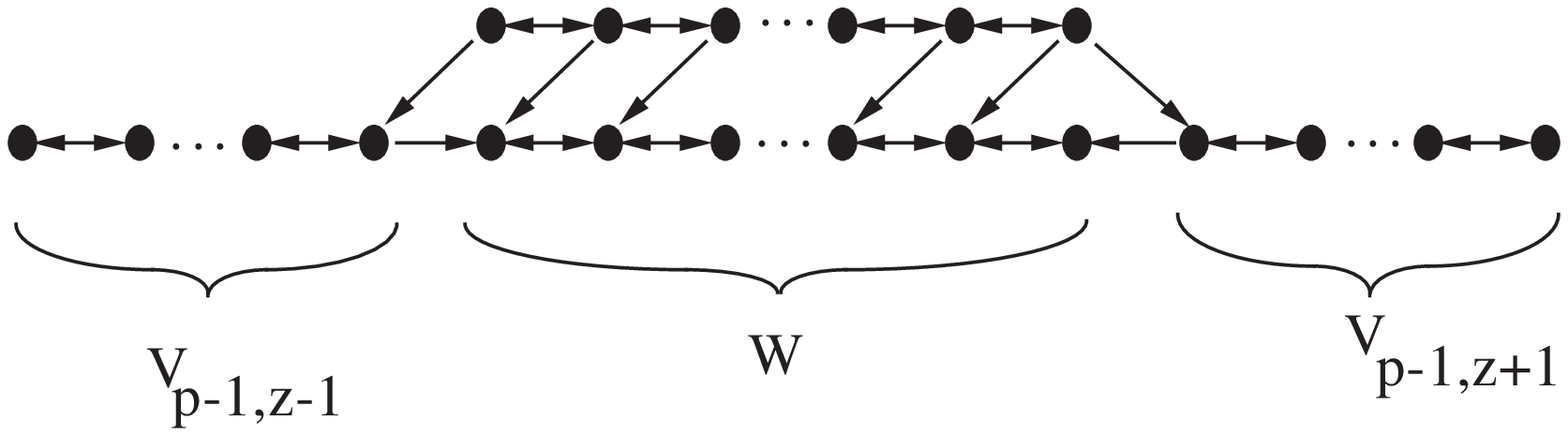}
   \vspace{-2.2in}
 \caption{Indecomposable \rep $I_z^p$}
\label{fig:indecomp}
\end{figure}

The irreps of $\U$ at roots of unity are well
known (see e.g. \cite{keller}, whose notations we largely follow), 
and we list some facts. 
Let 
\beq
q=e^{2 \pi i n/m}
\eq
for positive relatively prime integers $m,n$  
and define $M=m$ if $m$ is odd, and $M=m/2$ if $m$ is even. Then 
it is consistent and appropriate in our context to set 
\beq
(X^{\pm})^M = 0
\eq
(if one uses $q^{H}$ instead of $H$, then $(X^{\pm})^M $ is central). 
All finite dimensional irreps are highest weight (h.w.) representations 
with dimension $d\leq M$. There are two types of irreps:
\begin{itemize}
\item{$V_{d,z} = \{e_h^j;\quad j= (d-1) + \frac{m}{2n}z, \quad 
                   h=j, j-2, \dots, -(d-1) + \frac{m}{2n}z \}$ 
             with dimension $d$, for any 
       $1 \leq d \le M$ and $z \in \Z$, where $H e_h^j = h e_h^j$}
\item{$I^1_z$ with dimension $M$ and h.w. 
        $(M-1) +  \frac{m}{2n}z $, for
        $z \in \compl \setminus \{\Z+\frac{2n}{m}r, 1 \leq r \leq M-1\}$. }
\end{itemize}
Note that in the second type, $z \in \Z$ is allowed, in which case we 
will write $V_{M,z} \equiv I^1_z$ for convenience. We will concentrate
on the $V_{d,z}$ -- \reps from now on.
Furthermore, the fusion rules at roots of unity state that 
$V_{d,z} \tens V_{d',z'}$ decomposes into 
$\oplus_{d''} V_{d'',z+z'} \bigoplus_p I_{z+z'}^p$ where $I_z^p$ 
are the well--known reducible, but indecomposable 
representations of dimension $2M$, see  figure \ref{fig:indecomp}
and \cite{keller}. 
If $q$ is {\em not} a root of unity, then the universal 
${\cal R} \in {\cal U} \tens {\cal U}$ given by
\beq
{\cal R} = q^{\frac 12 H \tens H}\sum_{l=0}^{\infty} 
q^{-\frac 12 l(l+1)}\frac{(q-
q^{-1})^l}{[l]!} q^{l H/2}(X^+)^l \tens q^{-l H/2} (X^-)^l   \label{univ_R}
\eq
defines the quasitriangular structure of $\U$. It 
satisfies e.g. 
\beq
\sigma (\Del(u)) = {\cal R} \Del(u) {\cal R}^{-1}, \quad u \in {\cal U}
\eq
where $\sigma(a\tens b) = b \tens a$. 
We will only consider \reps with dimension $\leq M$; then ${\cal R}$ 
restricted to such \reps is well defined for roots of unity as well,
since the
sum in (\ref{univ_R}) only goes up to $(M-1)$. Furthermore
\beq
{\cal R}^{\ast} = ({\cal R})^{-1}.
\eq
To see this, (\ref{invol_help}) is useful.

Let us consider a hermitian invariant inner product  
$(u,v)$ for $u,v \in V_{d,z}$. A hermitian inner 
product satisfies $( u, \l v) = \l (u,v) = (\l^{c.c.} u, v)$ for 
$\l \in \compl$, $(u,v)^{c.c.} = (v,u)$, and it is invariant if
\beq
(u, x \cdot v) = (x^{\ast}\cdot u, v),
\eq
i.e. $x^{\ast}$ is the adjoint of $x$. If $(\ ,\ )$ is also positive
definite, we have a unitary representation.

\begin{prop}
The representations $V_{d,z}$ are 
unitarizable w.r.t.  $U_q(SO(2,1))$ if and only if 
\beq
(-1)^{z+1} \sin(2 \pi n k/m) \sin(2 \pi n(d-k)/m) >0
\eq
 for all  $k=1,...,(d-1)$. 

For $d-1 < \frac{m}{2n}$, this holds precisely if $z$ is odd. 
For $d-1 \geq \frac{m}{2n}$, 
it holds for isolated values of $d$ only, i.e. if it holds for $d$, 
then it (generally) does 
{\em not} hold for $d \pm1, d \pm 2,\dots$.

The representations $V_{d,z}$ are 
unitarizable w.r.t.  $U_q(SU(2))$ if
$z$ is even and $d-1 < \frac{m}{2n}$.
 \label{unitary_2d}
\end{prop}

\begin{proof} 
Let ${e^j_h}$  be a basis of $V_{d,z}$ with highest weight $j$.
After a straightforward calculation, invariance implies
\beq
\((X^-)^k \cdot e^j_j, (X^-)^k \cdot e^j_j\) = 
  (-1)^k [k]! [j] [j-1] ... [j-k+1] 
\(e^j_j, e^j_j\)
\eq
for $k=1, ... , (d-1)$, where $[n]! = [1] [2] ... [n]$.
Therefore we can have a positive definite inner product 
$(e^j_h, e^j_l) = \d_{h,l}$  if and only if 
$a_k \equiv (-1)^k [k]! [j] [j-1] ... [j-k+1]$ is a positive number for 
all $k=1,... ,(d-1)$, in which case 
$e^j_{j-2k} = (a_k)^{-1/2} (X^-)^k \cdot e^j_j$. 

Now $a_{k}= -[k] [j-k+1] a_{k-1}$, and
\berr
-[k] [j-k+1]  &=& -[k] [d-k+\frac{m}{2n}z] = -[k] [d-k] e^{i \pi z} \\                  
                     &=& (-1)^{z+1} \sin(2 \pi n k/m) \sin(2 \pi n(d-k)/m) 
\frac{1}{sin(2 \pi n/m)^2},
\err
since $z$ is an integer. Then the assertion follows.
The compact case is known \cite{keller}.
\end{proof}

In particular, all of them are finite dimensional, and clearly if $q$ is
not a root of unity, none of the \reps are unitarizable. 

We will be particularly interested in the case of (half)integer 
representations of type $V_{d,z}$ and $n=1, m$ even, for reasons to 
be discussed below. Then $d-1 < \frac{m}{2n} = M$ always holds, and 
{\em the $V_{d,z}$ are unitarizable if and only if $z$ is odd}. These 
representations are centered around $M z$, with dimension 
less than or equal to $M$.


Let us compare this with the classical case. For the 
Anti--de Sitter group $SO(2,1)$, $H$ is 
nothing but the energy (cp. section 3). At $q=1$, the unitary 
irreps of $SO(2,1)$ are lowest weight representations with 
lowest weight $j > 0$ resp. highest weight representations with 
highest weight $j<0$. For any given such lowest resp. highest weight 
we can now find a {\em finite dimensional} unitary 
representation with the same lowest resp. highest weight, 
provided $M$ is large enough (we only consider (half)integer $j$ here).  
These are unitary 
representations which for low energies look like the classical  one--
particle representations, but have an intrinsic high--energy cutoff
if $q \neq 1$, which goes to infinity as $q \rightarrow 1$.
The same will be true in the 4 dimensional case.

So far we only considered  what could be called one--particle 
representations. To talk about many--particle representations, 
there should be a tensor product of 2 or more such irreps, which gives a 
unitary representation as well and agrees with the classical case for
low energies. 

Since $\U$ is a Hopf  algebra, there is a natural notion of a tensor 
product of two representations, given by the coproduct $\Del$. 
However, it is not unitary a priori. 
As mentioned above, the tensor product of two irreps of type 
$V_{d,z}$ is
\beq
V_{d,z} \tens V_{d',z'} = \oplus_{d''} V_{d'',z+z'}  
\bigoplus_{p=r,r+2,\dots}^{d+d'-M} I_{z+z'}^p  
\label{fusion}
\eq
where $r=1$ if $d+d'-M$ is odd or else $r=2$, and
$I_z^p$ is a indecomposable representation of dimenson $2M$ 
whose structure is shown in figure \ref{fig:indecomp}.
The arrows indicate the 
rising and lowering operators. 

In the case of $U_q(SU(2))$, one usually defines a truncated tensor 
product $\obar{\tens}$ by omitting all indecomposable $I_z^p$ 
representations \cite{mack_schom}. Then the remaining \reps are 
unitary w.r.t. $U_q(SU(2))$; $\obar{\tens}$ is
associative only from the representation  theory point of view
\cite{mack_schom}.

This is not  the right thing to do for $U_q(SO(2,1))$. 
Let $n=1$ and $m$ even, and consider e.g. $V_{M-1,1}\tens V_{M-1,1}$. 
Both factors have lowest energy $H=2$, 
and the tensor product of the two corresponding {\em classical} \reps 
is the sum 
of \reps with lowest weights $4,6,8,\dots$ . 
In our case, these weights 
are in the $I_z^p$  representations, while the $V_{d'',z''} $ have 
$H \geq M \rightarrow \infty$ and are not unitarizable. 
So we have to keep the $I_z^p$'s 
and throw away the $V_{d'',z''} $'s in (\ref{fusion}). A priori however,
the $I_z^p$'s are 
not unitarizable, either. To get a unitary tensor product, note that 
as a vector space,
\beq	
I_z^p = V_{p-1, z-1} \oplus W \oplus V_{p-1,z+1} 
\label{I_decomp}
\eq
(for $p\neq 1$) where 
\beq
W=V_{M-p+1,z} \oplus V_{M-p+1,z}  \label{decomp}
\eq
as vector space. Now
$(X^+)^{p-1} \cdot v_l$ is a lowest weight vector where $v_l$  is the 
vector with lowest weight of $I_p^z$, 
and similarly $(X^-)^{p-1} \cdot v_h$ is a 
highest weight vector with $v_h$ being the 
vector of $I_p^z$ with highest weight (see figure \ref{fig:indecomp}). 
It is therefore consistent to consider the submodule of 
$I_p^z$ generated by $v_l$, and factor out its submodule generated by 
$(X^+)^{p-1} \cdot v_l$; the result is an irreducible \rep  
equivalent to $V_{p-1,z-1}$  realized on the left summand in 
(\ref{I_decomp}).
Similarly, one could consider the submodule of
$I_p^z$ generated by $v_h$, factor out its submodule generated by
$(X^-)^{p-1} \cdot v_h$, and obtain an irreducible \rep 
equivalent to $V_{p-1,z+1}$.
In short, one can just "omit" $W$ in (\ref{I_decomp}).
The two $V$ - type \reps obtained this way are unitarizable 
provided $n=1$ and $m$ 
is even, and one can either keep both 
(notice the similarity with band structures in solid--state 
physics), or for simplicity keep the 
low--energy part only, in view of the physical application we have in mind. 
We therefore  define a truncated tensor product  as
\begin{definition}
For $n=1$ and even $m$,
\beq 
V_{d,z} \ttens V_{d',z'} := \bigoplus_{\tilde{d}=r, r+2, \dots}^{d+d'-M} 
V_{\tilde{d},z+z'-1}
\label{fusion_t}
\eq
\end{definition}

This can be stated as follows: 
Notice that any representation 
naturally decomposes as a vector space into sums of $V_{d,z}$'s, 
cp. (\ref{decomp}); the 
definition of $\ttens$ simply means that only the smallest value of 
$z$ in this decomposition is kept, which is the submodule of irreps 
with lowest weights less than or equal to $\frac m{2n} (z+z'-1)$. 
(Incidentally, $z$ is the eigenvalue of $D_3$ in the {\em classical} 
$su(2)$ algebra generated by 
$\{D^{\pm} = \frac{(X^{\pm})^M}{[M]!}, 2D_3 = [D^+, 
D^-]\}$, where $\frac{(X^{\pm})^M}{[M]!}$ is understood by some 
limes procedure). 
With this in mind, it is obvious that $\ttens$ is associative: both in 
$(V_1\ttens V_2) \ttens V_3$ and in $V_1 \ttens (V_2 \ttens V_3)$, 
the result is simply the $V$'s with minimal $z$, which is the 
{\em same} space, because the ordinary tensor product is 
associative and $\Del$ is coassociative. This is in contrast 
with the "ordinary" truncated tensor product $\obar{\tens}$
\cite{mack_schom}. Of course, one could give a similar definition for
negative--energy representations. See also
definition \ref{trunc_tensor} in the case of $U_q(SO(2,3))$.

$V_{d,z} \ttens V_{d',z'}$ is unitarizable
if all the $V$ 's on the rhs of (\ref{fusion_t}) are unitarizable. 
This is certainly true if $n=1$ and $m$ is even. In all other cases,
there are no terms on the rhs of  (\ref{fusion_t}) if the factors
on the lhs are unitarizable, since no $I_z^p$ -- type \reps
are generated (they are too large). 
This is the reason why we concentrate on this case, 
and furthermore on
$z=z'=1$ which corresponds to low--energy representations.
Then $\ttens$ defines a two--particle Hilbert space with the 
correct classical limit. To summarize, we have the following:

\begin{prop}
$\ttens$ is associative, and 
$V_{d,1} \ttens V_{d',1}$ is unitarizable.
\end{prop}

How an inner product is induced from the single --
particle Hilbert spaces is explained in \cite{thesis}.

\section{The quantum group $U_q(SO(2,3))$ }

In order to generalize the above results to the 4 -- dimensional case, 
one has to use the general machinery of quantum groups, which 
is  briefly reviewed (cp. e.g. \cite{CH_P}):
Let $q \in \compl$ and $A_{ij}=2\frac{(\a_i, \a_j)}{(\a_j,\a_j)}$ be 
the Cartan matrix of a classical simple Lie algebra ${g}$ of rank $r$, 
where $( , )$ is the Killing form and $\{\a_i,\quad i=1,\dots,r \} $ 
are the simple roots. 
Then the {\em quantized universal enveloping algebra} 
$U_q({g})$ is the Hopf  algebra generated 
by the elements 
$\{ X^{\pm}_i, H_i; \quad i=1,\dots,r\}$ and relations 
\cite{FRT,jimbo,drinfeld}
\berr
\[H_i, H_j\]      &=& 0                       \nonumber  \\
\[H_i, X^{\pm}_j\] &=& \pm A_{ji} X^{\pm}_j,   \nonumber   \\
\[X^+_i, X^-_j\]   &=& \d_{i,j} \frac{q^{d_i H_i} 
                          -q^{-d_i H_i}}{q^{d_i}-q^{-d_i}}  
                    = \d_{i,j} [H_i]_{q_i},     \nonumber \\
\sum_{k=0}^{1-A_{ji}} &\ & 
           \[\begin{array}{c} 1-A_{ji} \\ k\end{array}\]_{q_i}(X^{\pm}_i)^k 
         X^{\pm}_j (X^{\pm}_i)^{1-A_{ji}-k} = 0, \quad i\neq j
\label{UEA}
\err
where $d_i = (\a_i, \a_i)/2, \quad q_i = q^{d_i}, \quad [n]_{q_i} 
           = \frac{q_i^n-q_i^{-n}}{q_i-q_i^{-1}}$ and 
\beq
\[ \begin{array}{c} n\\m \end{array} \]_{q_i} = 
             \frac{[n]_{q_i}!}{[m]_{q_i}! [n-m]_{q_i} !} .
\eq
The comultiplication is given by 
\berr
\Del(H_i)          &=& H_i \tens 1 + 1 \tens H_i , \nonumber \\
\Del(X^{\pm}_i) &=&  X^{\pm}_i \tens q^{d_i H_i/2} + q^{-d_i H_i/2} 
\tens  X^{\pm}_i. 
\err
Antipode and counit are  
\berr 
S(H_i)    &=& -H_i, \nonumber \\
S(X^+_i)  &=& -q^{d_i} X^+_i, \quad S(X^-_i)  = -q^{-d_i} X^-_i, \nonumber \\
\eps(H_i) &=& \eps(X^{\pm}_i)=0.
\err
(we use the conventions of \cite{kirill_resh}, which differ slightly from 
e.g. \cite{CH_P}.)

\begin{figure}
 \epsfxsize=5in
  \vspace{-2.5in}
   \hspace{3.5in}
  \epsfbox{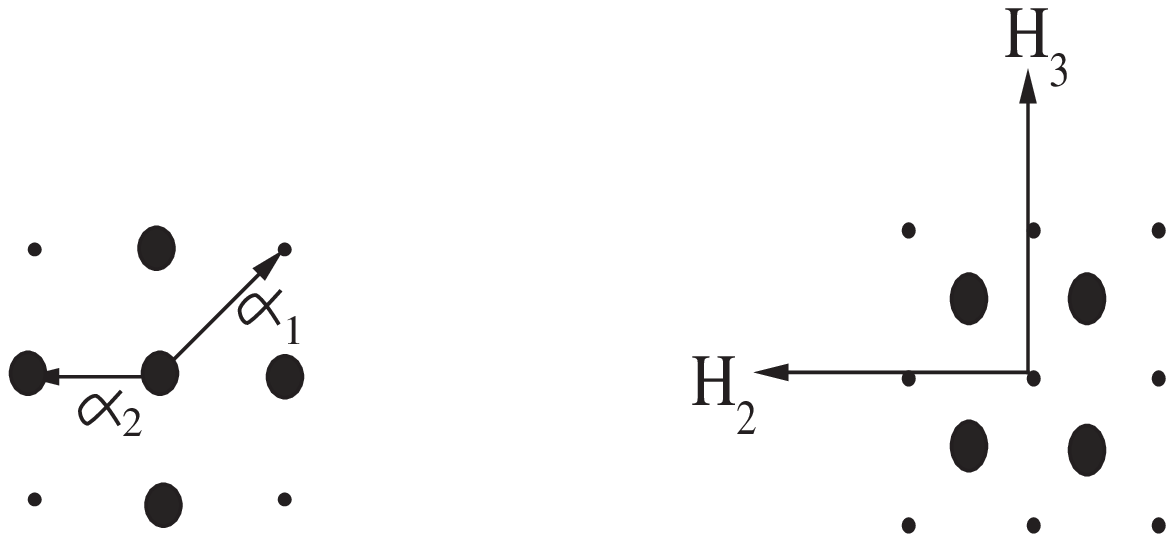}
   \vspace{-2.5in}
\caption{vector and spinor \reps of $SO(2,3)$}
\label{fig:so23reps}
\end{figure}

For $\U \equiv U_q(SO(5, \compl))$, $r=2$ and 
\beq
A_{ij} = \(\begin{array}{cc}    2 & -2 \\
                                -1 & 2 \end{array} \), 
         \qquad (\a_i, \a_j) = \( \begin{array}{cc} 2 & -1 \\
                                                   -1 &  1  \end{array} \) ,
\eq
so $d_1 = 1$, $d_2= 1/2$, to have the standard physics normalization 
(a rescaling of $(\ ,\ )$ can be absorbed by a redefinition of $q$). 
The weight diagrams of the vector 
and the spinor \reps are  given in figure \ref{fig:so23reps} for
illustration. The Weyl element  is 
$\rho = \frac 12 \sum_{\a >0} \a = \frac 32 \a_1 + 2 \a_2$.

The possible reality structures on $\U$ 
have been investigated in \cite{real_luk}. 
As in section 2,  in order to obtain finite dimensional unitary 
representations, $q$ must be a root of unity. Furthermore, on 
physical grounds we insist upon having positive--energy 
representations; already in the classical case, that rules out e.g. 
$SO(4,1)$, cp. the discussion in \cite{fronsdal}.
It appears that then there is only one possibility,
namely
\berr
(H_i)^{\ast} &=& H_i, \quad (X^+_1)^{\ast} = -X^-_1, \quad 
(X^+_2)^{\ast} = X^-_2,   	\label{star}  \\ 
(a\tens b)^{\ast} &=& b^{\ast} \tens a^{\ast},  \nonumber \\
(\Del(u))^{\ast} &=& \Del(u^{\ast}), \quad (S(u))^{\ast}= S(u^{\ast}),  
\label{star_Del}
\err
for $|q|=1$, which corresonds to the Anti--de Sitter group 
$U_q(SO(2,3))$.
Again with $E \equiv d_1 H_1+d_2 H_2$, 
$(-1)^E x^{\ast} (-1)^E = \theta(x^{c.c.})$ 
where $\theta$ is the usual Cartan--Weyl involution corresponding to 
$U_q(SO(5))$.

Although it will not be used it in the present paper, this algebra has the 
very important property of being {\em quasitriangular}, i.e. there 
exists a universal ${\cal R}$ $\in \U \tens \U$. It 
satisfies ${\cal R}^{\ast} = ({\cal R})^{-1}$,
which can be seen e.g. from uniqueness theorems, cp. 
\cite{tolstoi,CH_P}. 
In the mathematical literature, usually a rational version of the 
above algebra, i.e. using $q^{d_i H_i}$ instead of $H_i$ is considered. Since 
we are only interested in specific representations, 
we prefer to work with $H_i$. We essentially work in the "unrestricted"
specialization, i.e. the divided powers
$(X_i^{\pm})^{(k)} = \frac{(X_i^{\pm})^k}{[k]_{q_i}!}$ 
are not included if $[k]_{q_i}=0$, 
although our results will only concern \reps
which are small enough so that the distinction is not relevant.

Often the following generators are more useful: 
\beq
h_i = d_i H_i, \quad e_{\pm i} =  \sqrt{[d_i]}X^{\pm}_i,
\eq
so that 
\berr
\[h_i, e_{\pm j}\] &=& \pm (\a_i, \a_j) e_{\pm j}, \nonumber \\
\[e_i, e_{-j}\]       &=&  \d_{i,j} [h_i] .
\err
In the present case,
i.e. $h_1= H_1, h_2 = \frac 12 H_2$, 
$e_{\pm 1} =  X^{\pm}_1$ and 
$ e_{\pm 2} = \sqrt{[\frac 12]}X^{\pm}_2$.

So far we only have the generators corresponding to the simple roots. 
A  Cartan--Weyl basis corresponding to all roots can be obtained 
e.g. using the braid group action introduced by Lusztig \cite{lusztig},
(see also \cite{CH_P,jantzen}) resp. the quantum Weyl group 
\cite{kirill_resh,soib,lev_soib,CH_P}. 
If $\omega = \tau_{i_1} ... \tau_{i_N}$ is a reduced expression for 
the longest
element of the Weyl group where $\tau_i$ is the reflection along $\a_i$, then 
$\{\a_{i_1}, \tau_{i_1} \a_{i_2}, ..., \tau_{i_1} ... 
 \tau_{i_{N-1}} \a_{i_N}\}$
is an ordered set of positive roots. 
We will use $\omega = \tau_1 \tau_2 \tau_1 \tau_2$ and denote them
$\b_1 = \a_1, \b_2 = \a_2, \b_3 = \a_1 + \a_2, \b_4 = \a_1 + 2\a_2$.
A Cartan--Weyl basis of root vectors of $\U$ can then be defined as 
$\{e_{\pm 1}, e_{\pm 3}, e_{\pm 4}, e_{\pm 2 }\} = \{e_{\pm 1}, T_1 
e_{\pm 2}, T_1 T_2 e_{\pm 1}, T_1 T_2 T_1 e_{\pm 2}\}$ and similarly 
for the $h_i$ 's, where the $T_i$  represent the braid group on $\U$ 
\cite{lusztig}:
\berr
T_i(H_j)   &=& H_j - A_{ij} H_i, \quad T_i X_i^+ = -X_i^- q_i^{H_i},
                      \nonumber \\
T_i(X_j^+) 
         &=& \sum_{r=0}^{-A_{ji}}(-1)^{r-A_{ji}}q_i^{-r}(X_i^+)^{(-A_{ji}-r)} 
          X_j^+ (X_i^+)^{(r)}  \nonumber \\
\err
where $T_i(\theta(x^{c.c.})) = \theta(T_i(x))^{c.c.}$. We  find
\berr
e_3 &=& q^{-1} e_2 e_1 - e_1e_2, 
       \quad e_{-3} = q e_{-1}e_{-2} - e_{-2} e_{-1}, 
        \quad h_3 = h_1 + h_2 ,  \nonumber \\
e_4 &=& e_2 e_3 - e_3 e_2, 
       \quad e_{-4} = e_{-3} e_{-2} - e_{-2} e_{-3}, 
       \quad h_4 = h_1 + 2 h_2. \label{rootvect}
\err
Similarly one defines the root vectors 
$X^\pm_{\b_l}$.
This can be used to obtain a Poincar\'e--Birkhoff--Witt basis of 
$\U = \U^- \U^0 \U^+$ where $\U^{\pm}$ is generated by the 
$X^{\pm}_i$ and $\U^0$ by the $H_i$: 
for $\underline{k}:=(k_1,\dots, k_N)$ where $N$ is the number of 
positive roots, let 
$X^+_{\underline{k}} = X_{\b_1}^{+k_1} \dots X_{\b_N}^{+k_N}$. 
Then the $X^{\pm}_{\underline{k}}$ form
a P.B.W. basis of $\U^{+}$, and similarly for $\U^-$ \cite{lusztig_90} 
(assuming $q^4\neq 1$).

Up to a trivial automorphism, (\ref{rootvect}) agrees with the basis used in 
\cite{lukierski}. 
The identification of the usual generators of the 
Poincar\'e group has also been given there and will not be repeated here,
except for pointing out that $h_3$ is the energy and $h_2$ is a 
component of angular momentum, see also \cite{fronsdal}. 
All of the above form 
$U_{\tilde{q}}(SL(2,\compl))$ subalgebras with appropriate $\tilde{q}$ 
(but not as coalgebras), because the $T_i$'s are algebra 
homomorphisms. The reality structure is 
\beq
e_1^{\ast} = -e_{-1}, \quad e_2^{\ast}= e_{-2}, \quad 
e_3^{\ast} = -e_{-3}, \quad e_4^{\ast}= -e_{-4}.
\eq
So the set $\{e_{\pm 2},h_2\}$ generates a  
$U_{\tilde{q}}(SU(2))$ algebra, and the other three 
$\{e_{\pm\a}, h_{\a}\}$ generate noncompact $U_{\tilde{q}}(SO(2,1))$ 
algebras, as discussed in section 2. 


\section{Unitary representations of $U_q(SO(2,3))$ and 
$U_q(SO(5))$}

In this section, we consider representations of $U_q(SO(2,3))$ and 
show that
for suitable roots of unity $q$, the irreducible
positive resp. negative energy representations are again unitarizable, 
if the highest resp. lowest weight lies in some "bands" in weight space. 
Their structure for low energies is 
exactly as in the classical case including the appearance of "pure gauge" 
subspaces for spin bigger than or equal to 1 in the "massless" case, 
which have to be factored out to obtain the 
physical, unitary representations. At high energies, there is
an intrinsic cutoff.

From now on $q=e^{2\pi i n/m}$.
Most facts about representations of quantum groups 
we will use can be found e.g. in \cite{deconc_kac}.
It is useful to consider the Verma modules $M(\l)$ for a highest weight
$\l$, which is the (unique) $\U$ - module having a highest weight 
vector $w_{\l}$ such that
\beq
\U^+ w_{\l} = 0, \quad H_i w_{\l} = \frac{(\l,\a_i)}{d_i} w_{\l},
\label{weights}
\eq
and the vectors $X^-_{\underline{k}} w_{\l}$ form a P.B.W. basis of $M(\l)$. 
On a Verma module, one can define a unique invariant
inner product $(\ ,\ )$, which is hermitian and satisfies 
$(w_{\l}, w_{\l}) =1$ and  
$(u, x \cdot v) = (\theta(x^{c.c.})\cdot u, v)$
for $x \in \U$, as in section 2 \cite{deconc_kac}.
$\theta$ is again the (linear) Cartan -- Weyl involution 
corresponding to $U_q(SO(5))$. 

The irreducible highest weight \reps can be obtained  from the corresponding 
Verma module by factoring out all submodules in the Verma module.
All submodules are null spaces w.r.t. the above inner product,
i.e. they are orthogonal to any state in $M(\l)$. Therefore one can 
consistently 
factor them out, and obtain a hermitian inner product on the 
quotient space $L(\l)$, which is the unique irrep 
with highest weight $\l$. 
To see that they are null, let 
$w_{\mu} \in M(\l)$ be in some submodule, 
so $w_{\l} \notin \U^+ w_{\mu}$. 
Now  for  $v \in \U^- w_{\l}$, it follows 
$(w_{\mu}, v) \in (\U^+ w_{\mu}, w_{\l}) = 0$.

The following discussion until the paragraph before 
Definition \ref{cpct_def} is
technical and may be skipped upon first reading.
Let $Q = \sum \Z \a_i$ be the root lattice and $Q^+ = \sum \Z_+ \a_i$
where $\Z_+ =\{0,1,2,\dots\}$.
We will write 
\beq
\l \succ \mu \quad \mbox{if} \quad \l - \mu \in Q^+.
\eq
For $\eta \in Q$, denote (see \cite{deconc_kac})
\beq
\Par(\eta) := \{\underline{k} \in \Z_+^N; \quad \sum k_i \beta_i = \eta\}.
\eq 
Let $M(\l)_{\eta}$ be the weight space with weight $\l - \eta$ in $M(\l)$. 
Then its dimension is given by 
$|\Par(\eta)|$. If $M(\l)$ contains a highest weight vector with weight
$\sigma$, then the multiplicity of the weight space 
$\(M(\l) / M(\sigma) \)_{\eta}$ is given by 
$|\Par(\eta)| - |\Par(\eta +\sigma- \l)|$,
and so on. We will see how this allows to 
determine the structure, i.e. the characters of the irreducible highest weight
representations.

As usual, the character of a representation $V(\l)$ with maximal weight $\l$
is the function on weight space defined by
\beq
\ch (V(\l)) = e^{\l} \sum_{\eta \in Q^+} \dim{V(\l)_{\eta}} e^{-\eta},
\eq
where $e^{\l-\eta}(\mu) := e^{(\l-\eta,\mu)}$, and
$V(\l)_{\eta}$ is the weight space of $V(\l)$ at weight $\l-\eta$.
The characters of inequivalent highest weight irreps (which are 
finite dimensional at roots of unity) are linearly independent.
Furthermore, the characters of Verma modules are the same as in the classical
case \cite{jantzen_m,deconc_kac},
\beq
\ch (M(\l)) = e^{\l} \sum_{\eta\in Q^+} |\Par(\eta)| e^{-\eta}.
\label{char_verma}
\eq

In general, the structure of Verma modules is quite complicated, and
the proper technical tool to describe it is its {\em composition series}.
For a $\U$ --module $M$ with a maximal weight,
consider a sequence of submodules
$ \dots \subset W_2 \subset W_1 \subset W_0 = M$ such that
$W_k / W_{k+1}$ is irreducible, and thus $W_k / W_{k+1}\cong L(\mu_k)$
for some $\mu_k$. 
(If the series is finite, it is sometimes called
a Jordan--H\"older series. For roots of unity it is infinite,
but this is not a problem for our arguments.
$W_{k+1}$ can be constructed inductively by fixing a
maximal submodule of $W_k$, e.g. as the sum of all but one
highest weight submodules of $W_k$).
While the submodules $W_k$ may not be unique,
it is obvious that we always have
$\ch (M) = \sum \ch (W_k / W_{k+1}) = \sum \ch (L(\mu_k))$.
Since the characters of irreps
are linearly independent, this decomposition of $\ch (M)$
is unique, and so are the subquotients
$L(\mu_k)$. We will study the composition series of the Verma module
$M(\l)$, in order to determine
the structure of the corresponding irreducible highest weight representation.

Our main tool to achieve this is
a remarkable formula by De Concini and 
Kac for 
$\det(M(\l)_{\eta})$, the determinant of the inner product matrix of  
$M(\l)_{\eta}$. Before stating it, we point out its use for
determining irreps:

\begin{lemma}
Let $w_{\l}$ be the highest weight vector in an irreducible highest weight
\rep $L(\l)$ with invariant inner product. If $(w_{\l}, w_{\l}) \neq 0$, 
then $(\;,\;)$ is non--degenerate, i.e.
\beq         \label{irrep_lemma}
det(L(\l)_{\eta}) \neq 0
\eq
for every weight space with weight $\l - \eta$ in $L(\l)$.
\end{lemma}
\begin{proof}
Assume to the contrary that there is a vector $v_{\mu}$ which is
orthogonal to all vectors of the same weight,
and therefore to all vectors of any weight.
Because $L(\l)$ is irreducible, there exists an $u \in \U$ such that
$w_{\l} =  u\cdot v_{\mu}$. But then
$(w_{\l},w_{\l}) = (w_{\l}, u\cdot v_{\mu})
= (u^{\dagger}\cdot w_{\l},v_{\mu}) =0$, which is a contradiction.
\end{proof}

Now we state the result of De Concini and Kac \cite{deconc_kac}:
\beq
\det(M(\l)_{\eta}) = \prod_{\b \in R^+} \prod_{m_{\b} \in \N} 
                    \([m_{\b}]_{d_\b} \frac{q^{(\l+\r-m_{\b}\b/2,\b)} - 
                  q^{-(\l+\r-m_{\b}\b/2,\b)}}{q^{d_{\b}} - 
                 q^{-d_{\b}}} \)^{|\Par(\eta - m_{\b}\b)|}    
\label{deconc}
\eq
in a P.B.W. basis for arbitrary highest weight $\l$, 
where $R^+$ denotes the positive roots (cp. section 3), and
$d_{\b} =  (\b,\b)/2$.

To get some insight, notice first of all that due to  
$|\Par(\eta-m_{\b}\b)|$  in the 
exponent, the product is finite. Now for some positive root $\b$, let
$k_{\b}$ be the smallest integer such that 
$D(\l)_{k_{\b},\b} := \([k_{\b}]_{d_\b} \frac{q^{(\l+\r-k_{\b}\b/2,\b)} - 
q^{-(\l+\r-k_{\b}\b/2,\b)}}{q^{d_{\b}} - q^{-d_{\b}}} \) = 0$, 
and consider 
the weight space at weight $\l - k_{\b}\b$, i.e. $\eta_{\b} = k_{\b}\b$.
Then $|\Par(\eta_{\b} - k_{\b}\b)| =1$ and
$\det(M(\l)_{\eta_{\b}})$ is zero, so there is a highest weight
vector $w_{\b}$ with weight $\l-\eta_{\b}$ 
(assuming for now that there is no other null state with weight larger than
$(\l-\eta_{\b})$). 
It generates a submodule which is again a Verma module 
(because $\U$ does not have zero divisors \cite{deconc_kac}), with 
dimension $|\Par(\eta - k_{\b}\b)|$ at weight
$\l - \eta$. This is the origin of the exponent. 
However the submodules generated by the $\w_{\b_i}$ are in general not 
independent, i.e. they may contain common highest weight vectors, 
and other highest weight vectors besides these $w_{\b_i}$ might exist. 
Nevertheless, all the highest weights $\mu_k$
in the composition series of $M(\l)$ are precisely
obtained in this way. This "strong linkage principle" will be proven
below, adapting the arguments in \cite{jantzen_m} for the classical case. 
While it is not a new insight
for the quantum case either \cite{dobrev,anderson},
it seems that no explicit proof has been given at least
in the case of even roots
of unity, which is most interesting from our point of view, as we will see.

To make the structure more transparent, 
let $\N_{\b}^T$ be the set of positive integers $k$ with $[k]_{\b} = 0$,
and $\N_{\b}^R$ the positive integers $k$ such that
$(\l + \r -\frac k2 \b,\b) \in \frac m{2n}\Z$.
Then
\beq
D(\l)_{k,\b} = 0 \Leftrightarrow k \in \N_{\b}^T
                              \quad\mbox{or} \quad k\in \N_{\b}^R.
\eq
The second condition is
$k = 2 \frac{(\l+\rho, \b)}{(\b, \b)} + \frac m{2n} \frac 2{(\b,\b)} \Z$,
which means that
\beq
\l-k \b = \sigma_{\b,l}(\l)
\eq
 where
$\sigma_{\b,l}(\l)$
is the reflection of $\l$ by a plane perpendicular to $\b$ through
$-\rho+\frac m{4nd_{\b}} l \b$, for some integer $l$. 
For general $l$, $\sigma_{\b,l}(\l) \notin \l+Q$; 
but $k$ should be an integer, so it is natural to define the
{\em affine Weyl group} $\W_{\l}$ of reflections in weight space 
to be generated by those reflections $\sigma_{\b_i,l_i}$ 
in weight space which map $\l$
into $\l+Q$. For $q=e^{2\pi i n/m}$,
two such allowed reflection planes perpendicular to $\b_i$  will differ 
by multiples of $\frac 12 M_{(i)}\b_i$; here $M_{(i)}=m$ for $d_i=\frac 12$,
while for $d_i=1$, $M_{(i)}= m$ or $m/2$ if $m$ is odd or even, respectively. 
Thus $\W_{\l}$ is generated by all reflections by these planes.
Alternatively, it is generated by the usual Weyl group
with reflection center $-\rho$, and translations by $M_{(i)}\b_i$.

Now the {\em strong linkage principle} states the following:
\begin{prop}
$L(\mu)$ is a composition factor of the Verma module $M(\l)$
if and only if $\mu$ is {\em strongly linked} to $\l$,
i.e. if there is a descendant sequence of weights
related by the affine Weyl group as
\beq
\l \succ \l_i = \sigma_{\b_i,l_i}(\l)\succ \dots \succ
 \l_{kj...i} = \sigma_{\b_k,l_k}(\l_{j...i}) = \mu
\label{W_link}
\eq
\label{h_w_theorem}
\end{prop}

\begin{proof}
The main tool to show this is the
formula (\ref{deconc}). Consider the inner product matrix 
$M_{\uk,\uk'}:=(X^-_{\uk}w_{\l},X^-_{\uk'}w_{\l})$;
it is hermitian, since $q$ is a phase. One can define an analytic continuation 
of it as follows: for the same P.B.W. basis, 
let $B_{\uk,\uk'}(q,\l):=(X^-_{\uk}w_{\l},X^-_{\uk'}w_{\l})_b$ 
be the matrix of the invariant {\em bilinear} form 
defined as in \cite{deconc_kac}, which is manifestly analytic in $q$ and $\l$
(one considers $q$ as a formal variable and replaces 
$q\rightarrow q^{-1}$ in the first argument of $(\;,\;)_b$). 
Then (\ref{deconc}) holds for all $q\in\compl$ and
arbitrary complexified $\l$ \cite{deconc_kac}. For $|q|=1$ 
and real $\l$, 
$B_{\uk,\uk'}(q,\l) = M_{\uk,\uk'}$. Let $\l'=\l + h\rho$ and
$q' = q e^{i\pi h}$ for $h \in \compl$; then $B_{\uk,\uk'}(q',\l')$ is 
analytic in $h$, and hermitian for $h\in\reals$. 
Furthermore, one can identify the modules                     
$M(\l')$ for different $h$  via the P.B.W. basis. In this sense,
the action of $X_i^{\pm}$ is  analytic in $h$ (it only depends on the
commutation relations of the $X^{\pm}_{\b}$). 
Now it follows 
(see theorem 1.10 in \cite{kato}, chapter 2 on matrices which are analytic
in $h$ and normal for real $h$) that 
the eigenvalues $e_j$ of $B_{\uk,\uk'}(q',\l')$ are analytic in $h$, 
and there exist analytic projectors $P_{e_{j}}$
on the eigenspaces $V_{e_{j}}$ which span the entire
vectorspace (except possibly at isolated points where some eigenvalues
coincide; for $h\in\reals$ however, the generic eigenspaces
are orthogonal and therefore remain independent even at such points).
These projectors provide an analytic
basis of eigenvectors of $B_{\uk,\uk'}(q',\l')$.
Now let
\beq
V_k := \bigoplus_{e_{j} \propto h^k} V_{e_{j}},
\eq
i.e. the sum of the eigenspaces whose eigenvalues
$e_{j}$ have a zero of order $k$ (precisely) at $h=0$.
Of course, $(V_k,V_{k'})_b=0$ for $k \neq k'$. 
The $V_k$ span the entire
space, they have an analytic basis as discussed, and have the following
properties:

\begin{lemma}
\label{order_lemma}
\begin{itemize}

  \item [1)] $(v_k,v)_b = o(h^k)$ for
             $v_k\in V_k$ and {\em any} (analytic) $v \in M(\l')$.
  \item [2)] $X_i^{\pm} v_k ={\displaystyle \sum_{l\geq k}} a_l v_l +
             {\displaystyle \sum_{l =1}^k}  h^l b_l v_{k-l}$
            for $v_l\in V_l$ and  $a_l, b_l$ analytic. In particular at $h=0$,
            \beq
              M^k := \oplus_{n \geq k} V_n
            \eq
             is invariant.
\end{itemize}
\end{lemma}

\begin{proof}
\begin{itemize}
  \item [1)] Decomposing $v$ according to $\oplus_l V_l$, only the (analytic)
      component in $V_k$ contributes in $(v_k,v)_b$, with a factor $h^k$
      by the definition of $V_k$ ($o(h^k)$ means at least $k$ factors of $h$).
  \item [2)] Decompose $X_i^{\pm} v_k = \sum_{e_{j}} a_{e_{j}} v_{e_{j}}$
      with analytic coefficients $a_{e_{j}}$
      corresponding to the eigenvalue $e_{j}$.
      For any $v_{e_{j}}$ appearing on the rhs, consider
      $(v_{e_{j}}, X_i^{\pm} v_k)_b = a_{e_{j}} (v_{e_{j}},v_{e_{j}})_b =
      c\ a_{e_{j}} e_{j} $ (with $c \neq 0$ at $h=0$, since
      $v_{e_{j}}$ might not be normalized).
      But the lhs is $(X_i^{\mp} v_{e_{j}}, v_k)_b = o(h^k)$ as shown above.
      Therefore $a_{e_{j}} e_{j} =o(h^k)$, which implies 2).
\end{itemize} \end{proof}

In particular,
$M(\l)/_{M^1}$ is irreducible and nothing but $L(\l)$.
(The sequence of submodules
$ ... \subset M^2 \subset M^1 \subset M(\l)$ is similar to the Jantzen
filtration \cite{jantzen_m}.)

By the definition of $M^k$, we have
\beq
\mbox{ord} (\det(M(\l)_{\eta})) = \sum_{k\geq 1} \dim M^k_{\l-\eta}
\eq
where $M^k_{\l-\eta}$ is the weight space of $M^k$ at weight $\l - \eta$,
and $\mbox{ord} (\det(M(\l)_{\eta}))$ is the order of the
zero of $\det(M(\l)_{\eta})$
as a function of $h$, i.e. the maximal power of $h$ it contains.
Now from (\ref{deconc}),
it follows that
\berr
\sum_{k\geq 1} \ch (M^k) &=& e^{\l} \sum_{\eta\in Q^+} (\sum_{k\geq 1}
       \dim M^k_{\l-\eta}) e^{-\eta}  \nn\\
  &=& e^{\l} \sum_{\eta\in Q^+} \mbox{ord}(\det(M(\l)_{\eta})) e^{-\eta}\nn\\
  &=& \sum_{\b\in \R^+} (\sum_{n\in \N_{\b}^T} + \sum_{n\in \N_{\b}^R})
          e^{\l} \sum_{\eta\in Q^+} |\Par(\eta-n\b)| e^{-\eta} \nn\\
  &=& \sum_{\b\in \R^+} (\sum_{n\in \N_{\b}^T} + \sum_{n\in \N_{\b}^R})
          \ch (M(\l-n\b)) \nn \label{char_calc}\\
\err
where we used (\ref{char_verma}).

Now we can show (\ref{h_w_theorem}) inductively.
Both the left and the right side of
$(\ref{char_calc})$ can be decomposed into a sum of characters of
highest weight irreps, according to their composition series.
These characters are linearly independent.
Suppose that $L(\l-\eta)$ is a composition factor of $M(\l)$. Then the
corresponding character is contained in the lhs of
$(\ref{char_calc})$, since $M(\l)/_{M^1}$ is irreducible.
Therefore it is also contained in one of the $\ch (M(\l-n\b))$
on the rhs. Therefore $L(\l-\eta)$ is a composition factor
of one of these
$M(\l-n\b)$, and by the induction assumption we obtain that
$\mu \equiv \l-\eta$ is strongly linked to $\l$ as in (\ref{W_link}).

Conversely, assume that $\mu$ satisfies (\ref{W_link}). By the induction
assumption, there exists a $n\in \N_{\b}^T \cup \N_{\b}^R$ such that
$L(\mu)$ is a subquotient of $M(\l-n\b)$. Then (\ref{char_calc}) shows
that $L(\mu)$ is a subquotient of $M(\l)$.
\end{proof}

Obviously this applies to other quantum groups as well.
In particular, we recover the well--known fact that for
$q=e^{2\pi i n/m}$, all $(X^-_i)^{M_{(i)}} w_{\l}$ are
highest weight vectors, and zero in an irrep.

\begin{figure}
\epsfxsize=5in
  \vspace{-2in}
  \hspace{2.6in}
\epsfbox{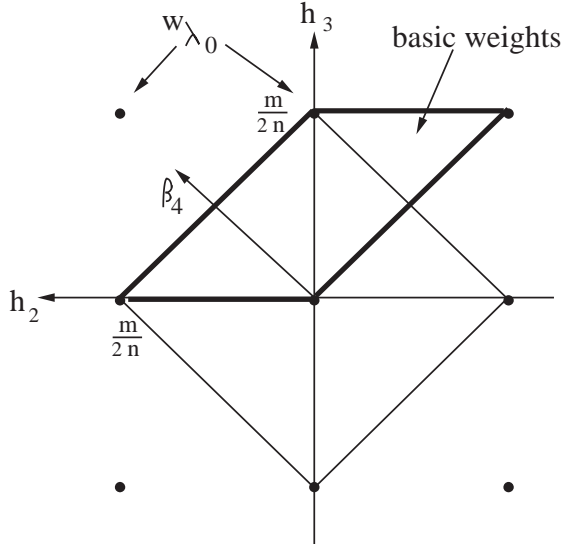}
  \vspace{-1.9in}
\caption{envelope of compact representations, basic weights and the 
lattice of $w_{\l_0}$}
\label{fig:basic}
\end{figure}

Now we can study the irreps of $U_q(SO(5))$ and 
$U_q(SO(2,3))$. First,
there exist remarkable nontrivial one dimensional 
\reps $w_{\l_0}$ with weights $\l_0 = \sum \frac m{2n} k_i \a_i$
for integers $k_i$. By tensoring
any \rep with $w_{\l_0}$, one obtains another \rep with the same structure, 
but all weights shifted by $\l_0$. We will see below that by such a shift, 
\reps which are 
unitarizable w.r.t. $U_q(SO(2,3))$ are in one to one correspondence 
with \reps which are unitarizable w.r.t. $U_q(SO(5))$.
It is therefore enough to consider highest weights in the following 
domain:

\begin{definition}
A weight $\l = E_0 \b_3 + s \b_2$ 
is called {\em basic} if 
\beq
0 \leq (\l, \b_3) = E_0 < \frac m{2n}, \quad 0 
  \leq (\l, \b_4) = (E_0+s) < \frac m{2n}.
  \label{basic}
\eq
In particular, $\l \succ 0$. It is said to be {\em compact} 
if in addition 
it is integral (i.e. $(\l, \b_i) \in \Z d_i$),  
\beq
s \geq 0 \quad \mbox{and} \quad (\l, \b_1) \geq 0.
  \label{cpct}
\eq
An irrep with compact highest weight will be called  compact.
\label{cpct_def}
\end{definition}
The region of basic weights is drawn in figure \ref{fig:basic},
together with the lattice of $w_{\l_0}$'s. The compact \reps are
centered around 0, and the (quantum) Weyl group acts on them
\cite{kirill_resh}, as classically (it is easy to see that the action of
the quantum Weyl group resp. braid group on the compact \reps 
is well defined at roots of unity as well).

A \rep with basic highest weight can be unitarizable w.r.t. 
$U_q(SO(5))$ (with conjugation $\theta^{c.c.}$)
only if all the $U_q(SU(2))$'s are unitarizable. For compact 
$\l$, all the $U_q(SU(2))$'s are indeed unitarizable according to section 2,
where $M_{(2)} = M_{(3)} = m$ and $M_{(1)} = M_{(4)} = m$ or $m/2$
if $m$ is odd or even, respectively. This alone
however is not enough to show that they are
unitarizable w.r.t. to the full group.

Although it may not be expected, there exist unitary 
representations with non--integral basic highest weight, namely 
for 
\beq
\l=\frac{m-1}2 \b_3 \quad \mbox{and } \quad 
\l=(\frac m2 -1)\b_3 + \frac 12 \b_2,
\label{nonint_cpct}
\eq
if $n=1$ and $m$ even. It follows from Proposition \ref{h_w_theorem} 
that they contain a highest weight vector at $\l-2\b_3$ and $\l-\b_3$
respectively, and all the
multiplicities in the irreps turn out to be one. 
Furthermore all $U_q(SU(2))$ modules in
$\b_1, \b_4$ direction have maximal length $M_{(1)} = m/2$,
which implies that they are unitarizable. 
The structure is that of shifted Dirac singletons which were 
already studied in \cite{dobrev}, and we will return to this later.

It appears that all other irreps must have integral highest weight 
in order to be unitarizable w.r.t. $U_q(SO(5))$.
If the highest weight 
is not compact, some of the $U_q(SU(2))$'s will not be unitarizable.
On the other hand, all irreps with compact highest weight are indeed
unitarizable: 

\begin{theorem}
The structure of the irreps $V(\l)$ with compact highest weight $\l$
is the same as in the classical case except if
\begin{itemize}
  \item [a)] $\l = (m/2-1-s)\b_3 + s \b_2$ for $s\geq 1$ and $\frac m{2n}$ 
        integer, where one additional highest weight state at weight 
       $\l - \b_4$ appears and no others, and 
  \item [b)] $\l=\frac{m-1}2 \b_3$ and 
       $\l=(\frac m2 -1)\b_3 + \frac 12 \b_2$  for $n=1$ and $m$ odd, 
       where one additional highest weight state at weight $\l-2\b_3$ 
        resp. $\l-\b_3$ appears and no others,
\end{itemize}
which are factored out in the irrep.
They are unitarizable w.r.t. $U_q(SO(5))$ 
(with conjugation $\theta^{c.c.}$).

The irreps with nonintegral highest weights 
(\ref{nonint_cpct}) discussed above are unitarizable as well.
\label{cpct_thm}
\end{theorem} 

\begin{proof}
The statements on the structure follow easily from Proposition
\ref{h_w_theorem}.

To show that these irreps are unitarizable, consider the
compact representation with highest weight $\l$ before factoring out 
the additional highest weight state,
so that the space is the same as classically. 
For $q=1$, they are known to be unitarizable, so the 
inner product 
is positive definite. Consider the eigenvalues of the (hermitian)
inner product matrix $M_{\uk,\uk'}$
as $q$ goes from 1 to $e^{2 \pi i n/m}$ along the unit circle.
The only way an eigenvalue could become negative is that it is zero
for some $q'\neq q$. This can only happen
if $q'$ is a root of unity, $q'=e^{2 i \pi n'/m'}$ with  $n'/m' < n/m$.
 then the "non--classical" reflection planes of $\W_{\l}$ 
are further away 
from the  origin and are relevant only in the case
$\l=\frac{m-1}2 \b_3$ for $n=1$ and $m$ odd;
but as pointed out above, no additional eigenvector 
appears in this case for $q' \neq q$.

Thus the eigenvalues might only become zero at $q$.
This happens precisely if a new highest weight vector appears,
i.e. in the cases listed. 
Since there is no null vector in the remaining irrep, all its
eigenvalues are positive by continuity.
\end{proof}

So far, all results were stated for highest weight modules; 
of course the 
analogous statements for lowest weight modules are true as well.

Now we want to find the "physical",  positive--energy \reps 
which are unitarizable w.r.t. $U_q(SO(2,3))$.
They are most naturally considered as lowest  weight representations, 
and can be obtained from the compact case by a shift,
as indicated above: if $V(\l)$ is a compact highest weight
representation, then 
\beq 
V(\l) \cdot \w\ := V(\l) \tens \w
\eq
with $\w \equiv w_{\l_0}, \l_0 = \frac m{2n} \b_3$ 
has lowest weight $\mu =  -\l + \l_0 \equiv E_0 \b_3 -s \b_2$ 
(short: $\mu = (E_0, s)$).
It is a positive--energy representation, i.e. the eigenvalues of 
$h_3$ are positive. 

\begin{figure}
 \epsfxsize=5in
   \vspace{-2in}
   \hspace{0.5in}
\epsfbox{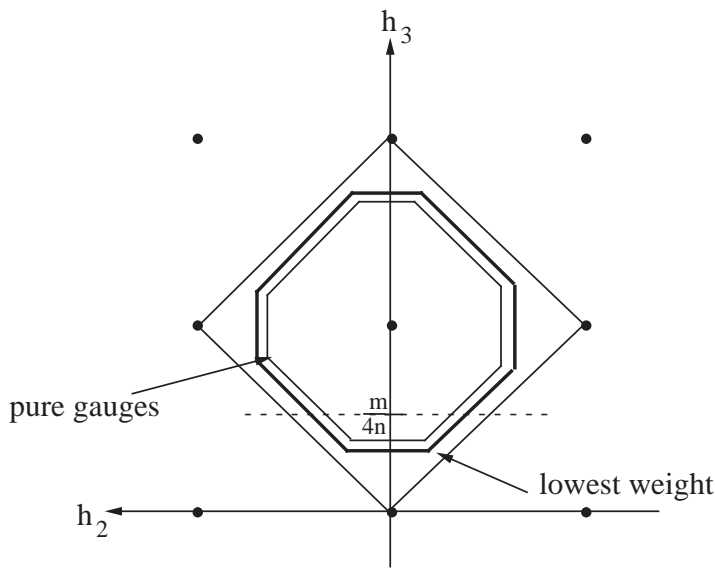}
   \vspace{-2.1in}
\caption{physical \rep with subspace of pure gauges 
  (for $\frac m{2n}$ integer), schematically. 
          For $h_3 \leq \frac m{4n}$,    
          the structure is the same as for $q=1$.}
\label{fig:gauge}
\end{figure}

For $\frac m{2n}$ integer, these \reps will correspond precisely to 
classical positive--energy representations
with the same lowest weight \cite{fronsdal}. The states with smallest 
energy $h_3$ corresponds to the particle at rest, so $E_0$ is the rest energy
and $s$ the spin. For $h_3 \leq m/4n$, the structure  is the same 
as classically, see figure \ref{fig:gauge}.
The irreps with nonintegral highest weights 
(\ref{nonint_cpct}) upon this shift correspond to the Dirac singleton 
\reps "Rac" with lowest weight $\mu=(1/2,0)$ and "Di" with $\mu=(1,1/2)$,
as discussed in \cite{dobrev}. 

If $\frac m{2n}$ is not integer, the weights of shifted compact \reps
are not integral. For $n=1$ and $m$ odd, the irreps in 
b) of theorem (\ref{cpct_thm}) now correspond to the singletons, 
again in argeement with \cite{dobrev}. We will see however that this case
does not lead  to an interesting tensor product. 

The case
$\mu = (s+1, s)$ for $s \geq 1$ and $\frac m{2n}$ integer will be 
called "massless"
for two reasons. First, $E_0$ is the smallest possible rest energy for 
a unitarizable \rep with given $s$ (see below). The main reason
however is the fact that as in the classical case
\cite{fronsdal}, an additional 
lowest  weight state with $E_0'=E_0+1$ and $s' = s-1$ appears, 
which generates a null subspace of what should be called 
"pure gauge" states. 
This corresponds precisely to the classical phenomenon in gauge theories,
which ensures that the massless photon, graviton etc. have only their
appropriate number of degrees of freedom (generally, the concept of mass
in Anti--de Sitter space is not as clear as in flat space.  Also notice 
that while
"at rest" there are actually still $2s+1$ states, the \rep 
is nevertheless reduced
by one irrep of spin $s-1$). 
In the present case,
all these representations are finite--dimensional! 

Thus we are led to the following.

\begin{definition}
An irreducible representation $V_{(\mu)}$ with lowest weight 
$\mu =  (E_0,s) \equiv E_0 \b_3 - s \b_2$ (resp. $\mu$ itself)
is called {\em physical} if it is unitarizable w.r.t. $U_q(SO(2,3))$
(with conjugation as in (\ref{star_Del})).

It is called {\em massless} if 
$E_0=s+1$  for $s \geq 1$, $s \in \frac 12 \Z$ and $\frac m{2n}$ integer.

For $n=1$, $V_{(\mu)}$ is called {\em Di} if $\mu=(1,1/2)$ and {\em Rac} if
$\mu=(1/2,0)$.

\label{physical}
\end{definition}


\begin{theorem}
The irreducible \rep $V_{(\mu)}$ with lowest weight $\mu$ is physical,
i.e. unitarizable w.r.t. $U_q(SO(2,3))$, 
if and only if the (shifted) irreducible represesentation 
with lowest weight
$\mu - \frac m{2n} \b_3$ is unitarizable w.r.t. $U_q(SO(5))$.

All $V_{(\mu)}$ where
$-(\mu - \frac m{2n} \b_3)$ is compact
are physical, as well as
the singletons Di and Rac. For $h_3 \leq \frac m{4n}$,
$V_{(\mu)}$ is obtained by factoring out from a (lowest weight)
Verma module a 
submodule with lowest weight $(E_0, -(s+1))$, except for the massless 
case, where one additional lowest weight state 
with weight $(E_0+1, s-1)$ appears, and for the Di and Rac, where one 
additional lowest weight state with weight $(E_0+1,s)$ and $(E_0+2,s)$ 
appears, respectively.
This is the same as for $q=1$, see figure \ref{fig:gauge}.
\label{noncpct_thm}
\end{theorem} 
For the singletons, this was already shown in \cite{dobrev}.

\begin{proof}
As mentioned before, we can write every vector in such a \rep uniquely 
as $a \cdot \w$, where $a$ belongs to a unitarizable
irrep of $U_q(SO(5))$. Consider the inner product
\beq
\<a\cdot \w, b\cdot\w \> \equiv (a,b),
\eq
where $(a,b)$ is the hermitian inner product on the {\em compact} (shifted) 
representation. Then
\berr
\<a\cdot\w, e_1 (b\cdot\w)\> &=& \<a\cdot\w,(e_1 \tens q^{h_1/2} + q^{-h_1/2}      
                    \tens e_1)b\tens\w \>  \nonumber \\
   &=& q^{h_1/2}|_{\w} (a, e_1 b) = i (a, e_1 b) 
\err
using $h_1|_{\w} = \frac m{2n}$. Similarly, 
\berr
\< e_{-1} (a\cdot\w), b\cdot\w\> &=& \< (e_{-1} \tens q^{h_1/2} +q^{-h_1/2} 
              \tens e_{-1})a\tens\w, b\tens\w\> \nonumber \\
   &=&  q^{-h_1/2}|_{\w} (e_{-1}a, b) = -i (e_{-1} a, b)
\err
because $\<,\>$ is antilinar in the first argument and linear in the second.
Therefore
\beq
\<a\cdot\w, e_1 (b\cdot\w)\>  = -\< e_{-1} (a\cdot\w), b\cdot\w\>.
\eq
Similarly $\<a\cdot\w, e_2 (b\cdot\w)\>  = \< e_{-2} (a\cdot\w), b\cdot\w\>$.
This shows that $x^{\ast}$ is indeed the adjoint of $x$ w.r.t. 
$\<\ ,\ \>$ which is positive definite,
because $(\ ,\ )$ is positive definite by definition. 
Theorem (\ref{cpct_thm}) now completes the proof.

\end{proof}

As a consistency check, one can see again from section 2 that all the 
$U_{\tilde{q}}(SO(2,1))$ resp. $U_{\tilde{q}}(SU(2))$ subgroups 
are unitarizable in these representations, but this 
is not enough to show the unitarizability for the full group. 
Note that for $n=1$, 
one obtains the classical one--particle \reps for given $s, E_0$ as 
$m \rightarrow \infty$.
We have therefore also proved the unitarizability at $q=1$
for (half)integer spin, which appears to be non trivial in itself 
\cite{fronsdal}. Furthermore, {\em all \reps obtained 
from the above by shifting
$E_0$ or $s$ by a multiple of $\frac mn$ are unitarizable as well}. 
One obtains in weight space 
a cell--like structure of representations which  are unitarizable w.r.t. 
$U_q(SO(2,3))$ resp. $U_q(SO(5))$. 

Finally we want to consider many--particle representations, i.e. find  a
tensor product such that the tensor product of unitary \reps 
is unitarizable, 
as in section 2. The idea is the same as there, the tensor product of 
2 such \reps
will be a direct sum of representations, of which we only keep the appropriate  
physical lowest weight "submodules".
To make this more precise, consider two physical irreps 
$V_{(\mu)}$ and $V_{(\mu')}$ as in Definition \ref{physical}. 
For a basis $\{u_{\l'}\}$ of lowest 
weight vectors in $V_{(\mu)} \tens V_{(\mu')}$ with physical $\l'$, 
consider the linear 
span $\oplus \U u_{\l'}$ of its lowest weight submodules, and 
let $Q_{\mu,\mu'}$ be the quotient of this after factoring out 
all proper submodules of the $\U u_{\l'}$. 
Let $\{u_{\l''}\}$ be a basis of lowest 
weight vectors of $Q_{\mu,\mu'}$.
Then $Q_{\mu,\mu'} = \oplus V_{(\l'')}$ where $V_{(\l'')}$ are the
corresponding (physical) irreducible lowest weight modules, 
i.e. $Q_{\mu,\mu'}$ is completely reducible.
Now we define the following:

\begin{definition} 
\label{trunc_tensor}
In the above situation, let $\{u_{\l''}\}$ be a basis of physical
lowest weight states of $Q_{\mu,\mu'}$, and let
$V_{(\l'')}$ be the corresponding physical lowest weight irreps. 
Then define  
\beq
V_{(\mu)} \ttens V_{(\mu')} := \bigoplus_{\l''} V_{(\l'')}
\eq
\end{definition} 
Notice that if $\frac m{2n}$ is not an integer, 
then the physical states
have non--integral weights, and  the full tensor product of two
physical irreps $V_{(\mu)} \tens V_{(\mu')}$ does not contain any 
physical lowest weights. Therefore 
$V_{(\mu)} \ttens V_{(\mu')}$ is zero in that case.

Again as in section 2, one might also include a second "band" 
of high--energy states. 

\begin{theorem} 
If all weights in the factors are integral, then
$\ttens$ is associative, and 
$V_{(\mu)} \ttens V_{(\mu')}$ is unitarizable w.r.t. $U_q(SO(2,3))$.
\end{theorem}
\begin{proof}
First, notice that the $\l''$ are all integral and none of them gives rise
to a  massless \rep or a singleton. Thus none of the $\U u_{\l'}$ 
contain a physical lowest weight vector
according to Proposition \ref{h_w_theorem}.
Also, lowest weight vectors for generic $q$ cannot disappear 
at roots of unity. 
Therefore $Q_{\mu,\mu'}$ contains all the physical 
lowest weight vectors of the
full tensor product. Furthermore, no physical lowest weight vectors
are contained in products of the form 
(discarded vectors)$\tens$ (any vectors).
Associativity now follows from the associativity of the full tensor 
product and the coassociativity of the coproduct, 
and the stucture 
for energies $h_3 \leq \frac m{4n}$ is the same as classically
(observe that since there are no massless representations, classically 
inequivalent physical \reps cannot recombine into 
indecomposable ones).
\end{proof}

In particular, none of the low--energy
states have been discarded. Therefore our definition is physically sensible,
and the case of $q=e^{2\pi i/m}$ with $m$ even appears to be most 
interesting physically.

\section{Conclusion}

We have shown that in contrast to the classical case, there exist
unitary \reps of noncompact quantum groups at roots of unity.
In particular, we have found finite dimensional unitary \reps 
of $U_q(SO(2,3))$ corresponding to all classical "physical"
representations, with the same structure 
at low energies as in the 
classical case. Thus they could be used to describe 
elementary particles with arbitrary spin.
This generalizes earlier results of \cite{dobrev} on the singletons.
Representations for many non--identical particles are found.

Apart from purely mathematical interest, this is very
encouraging for applications in QFT. In particular the appearance 
of pure gauge states should be a good guideline to construct
gauge theories on quantum Anti--de Sitter space. 
If this is possible, one should expect it to be finite in light of these
results.
However to achieve that goal,
more ingredients are needed, such as implementing a symmetrization axiom
(cp. \cite{schupp}),
a dynamical principle (which would presumably involve integration over
such a quantum space, cp. \cite{integral}), and efficient methods to
do calculations in such a context. These are areas of current research.

\section{Acknowledgements}
It is a pleasure to thank Bruno Zumino for many useful discussions, 
encouragement and support. I would also like to thank  
Chong--Sun Chu, Pei--Ming Ho and  Bogdan Morariu for useful comments, 
as well as  Paolo Aschieri, Nicolai Reshetikhin,
Michael Schlieker and Peter Schupp.
In particular,  I am grateful to V. K. Dobrev for pointing out the papers
\cite{dobrev} to me.
This work was supported in part by the Director, Office of
Energy Research, Office of High Energy and Nuclear Physics, Division of
High Energy Physics of the U.S. Department of Energy under Contract
DE-AC03-76SF00098 and in part by the National Science Foundation under
grant PHY-9514797.

\end{document}